\documentclass[twoside,12pt]{article}
\title{Rings whose associated  extended zero-divisor graphs are complemented  }

\date{}
\author{}
\usepackage[all]{xy}
\usepackage{amssymb,amsmath,latexsym}
\usepackage{theorem}
\usepackage{amsfonts}
\usepackage{ulem}
\usepackage{amscd}
\usepackage{amsxtra}
\usepackage{graphicx}
\usepackage{bm}
\usepackage{hyperref}
\usepackage{graphics} %inclusion de figures
\usepackage{graphicx} %inclusion de figures
\usepackage{pstricks,pst-node} %Graphiques
\usepackage{tikz} %Tikz !

\usepackage{tabularx} %Tableaux
\usepackage{multirow} %Gestion des lignes
\usepackage{multicol} %Gestion des colones
\usepackage{arydshln} %Lignes en pointillés
\usepackage{fancybox} %Boites
\usepackage{multicol} %Colonnes
\usepackage{array} %Tableaux maths
\usepackage{fancybox}

\usepackage{mathrsfs}
\usepackage{array,multirow,makecell}
\setcellgapes{1pt}
\makegapedcells

\newcolumntype{R}[1]{>{\raggedleft\arraybackslash }b{#1}}
\newcolumntype{L}[1]{>{\raggedright\arraybackslash }b{#1}}
\newcolumntype{C}[1]{>{\centering\arraybackslash }b{#1}}

\usetikzlibrary{arrows} 

\newcommand{\field}[1]{\mathbb{#1}}

\newcommand{\R}{\field{R}}

\newcommand{\Z }{\field{Z}}
\newcommand{\N }{\field{N}}

\theoremstyle{}\newtheorem{thm}{\bf Theorem}[section]
\theoremstyle{}\newtheorem{cor}[thm]{\bf Corollary}
\theoremstyle{}\newtheorem{lem}[thm]{\bf Lemma}
\theoremstyle{}
\theoremstyle{}
\theoremstyle{}\newtheorem{pro}[thm]{\bf Proposition}
\theoremstyle{}\newtheorem{exm}[thm]{\bf Example}
\theoremstyle{}
\theoremstyle{}\newcommand{\cqfd}{\hfill$\square$}

\def\pr{{\parindent0pt {\bf Proof.\ }}}

\def\cqfd
{\hspace{1cm}
	\rule{2mm}{2mm}%
	\medbreak%
	\par%
}

\def\ann{{\rm Ann}}
\def\nil{{\rm Nil}}

\begin{document}
	
	%%%%%%%%%%%%%%%%%%%%%%%%%%%%%%%%%%%%%%%%%%%%%%%%%%%%%%%%%
	%%%%%%%%%%%%%%%%%%%%%%%%%%%%%%%%%%%%%%%%%%%%%%%%%%%%%%%%%
	
	%%%%%%%%%%%%%%%%%%%%%%%%%%%%%%%%%%%%%%%%%%%%%%%%%%%%%%%%%
	%%%TITLE%%%%%%%%%%%%%%%%%%%%%%%%%%%%%%%%%%%%%%%%%%%%%%%%%
	\maketitle \vspace*{-1.5cm}
	
	%%%%%%%%%%%%%%%%%%%%%%%%%%%%%%%%%%%%%%%%%%%%%%%%%%%%%%%%%
	%%%%%%%%%%%%%%%%%%%%%%%%%%%%%%%%%%%%%%%%%%%%%%%%%%%%%%%%%
	%%%%%%%%%%%%%%%%%%%%%%%%%%%%%%%%%%%%%%%%%%%%%%%%%%%%%%%%%
	%%%NAMES%%%%%%%%%%%%%%%%%%%%%%%%%%%%%%%%%%%%%%%%%%%%%%%%%
	
	\begin{center}
		{\large\bf   Driss Bennis$^{1, a}$, Brahim El Alaoui$^{1, b}$, and  Raja L'hamri$^{1,e}$}
		
		\bigskip
		
		%%%%%%%%%%%%%%%%%%%%%%%%%%%%%%%%%%%%%%%%%%%%%%%%%%%%%%%%%
		%%%ADDRESSES%%%%%%%%%%%%%%%%%%%%%%%%%%%%%%%%%%%%%%%%%%%%%
		
		$^1$   Faculty of Sciences,  Mohammed V University in Rabat,  Morocco.\\
		\noindent    $^a$\,driss.bennis@um5.ac.ma; driss$\_$bennis@hotmail.com;  \\ $^b$\,brahimelalaoui0019@gmail.com; brahim$\_$elalaoui2@um5.ac.ma; \\   $^e$\,raja.lhamri@um5s.net.ma; rajaaalhamri@gmail.com  \\[0.2cm]  
	\end{center}
	\bigskip \bigskip
	%%%%%%%%%%%%%%%%%%%%%%%%%%%%%%%%%%%%%%%%%%%%%%%%%%%%%%%%%
	%%%%%%%%%%%%%%%%%%%%%%%%%%%%%%%%%%%%%%%%%%%%%%%%%%%%%%%%%
	%%%%%%%%%%%%%%%%%%%%%%%%%%%%%%%%%%%%%%%%%%%%%%%%%%%%%%%%%
	%%%ABSTRACT%%%%%%%%%%%%%%%%%%%%%%%%%%%%%%%%%%%%%%%%%%%%%%
	
	\noindent{\large\bf Abstract.}

	Let $R$ be a commutative ring with identity $1\neq 0$.
	 In this paper,  we continue the study started in  \cite{DJF} concerning  when the extended zero-divisor graph of $R$, $\overline{\Gamma}(R)$, is complemented. We also study when $\overline{\Gamma}(R)$  is uniquely complemented. We give a complete characterization of when $\overline{\Gamma}(R)$ of a finite ring is complemented. Various examples are given using the direct product of rings and idealizations of modules.\\
	

	
	%%%%%%%%%%%%%%%%%%%%%%%%%%%%%%%%%%%%%%%%%%%%%%%%%%%%%%%%%
	%%%%%%%%%%%%%%%%%%%%%%%%%%%%%%%%%%%%%%%%%%%%%%%%%%%%%%%%%
	%%%%%%%%%%%%%%%%%%%%%%%%%%%%%%%%%%%%%%%%%%%%%%%%%%%%%%%%%
	%%%Section 1 
	%%%%%%%%%%%%%%%%%%%%%%%%%%%%%%%%%%%%%%%%%%%%%%%%%%%%%%%%%
	%%%%%%%%%%%%%%%%%%%%%%%%%%%%%%%%%%%%%%%%%%%%%%%%%%%%%%%%%
	%%%%%%%%%%%%%%%%%%%%%%%%%%%%%%%%%%%%%%%%%%%%%%%%%%%%%%%%%
	
	\section{Introduction}
	Throughout the paper, $R$  will be a commutative ring  with  identity and $Z(R)$ be  its set of zero-divisors. Let  $x$ be an element  of  $R$, the annihilator of $x$  is defined as  $\ann_R(x):=\{y\in R/ \ xy=0\}$ and if there is no confusion,  we denote it simply by $\ann(x)$. For an ideal $I$ of $R$, $\sqrt{I}$ means the radical of $I$.  An element $x$ of $R$ is called nilpotent if $x^n=0$ for some positive integers $n$ and we denote $n_x$ its index of nilpotency; that is  the smallest integer $n$ such that $x^n=0$. The set of all nilpotent elements is denoted   $\nil(R):=\sqrt{0}$. The ring $\Z/n\Z$ of the residues modulo an integer $n$ will be denoted by $\Z_n$.  
	For a subset $X$ of $R$, we denote by  $X^*$ the set $X\setminus \{0\}$. \\
	Recall that the zero-divisor graph, denoted by  $\Gamma(R)$,  is the simple graph whose vertex set is the set of nonzero zero divisors, denoted by $Z(R)^*=Z(R)- \{0\}$,   and two distinct  vertices $x$ and $y$  are adjacent if and only if  $xy=0$.
	The extended zero-divisor graph, denoted by	 $\overline{\Gamma}(R)$,  is the simple graph who has the same vertex set like  $\Gamma(R)$ and two distinct  vertices $x$ and $y$ are adjacent if and only if $x^ny^m=0$ with $x^n\neq 0$ and $y^m\neq 0$  for some integer $n,m\in \N^*$.  We assume the reader has at least a basic familiarity with  zero-divisor graph theory.  For general background   on zero-divisor graph theory, we refer the reader to \cite{ ADL,DM,BAD,AM07,Axt,BI, DJF, BEFFR}.\\

This paper deals with complementeness and uniquely complementeness notions of graphs. A graph $G=(V,E)$ is said to be complemented if every vertex has an orthogonal. Namely, for every vertex $v$ in $G$ there exists a vertex $u$ in $G$ such that the edge $v-u$ is not a part of a triangle. In other words for every $v\in V$, there exists $u\in V$ such that $v$ is orthogonal to $u$, denoted by $v\perp u$.  $G$ is said to be uniquely complemented if it is complemented  and for any three vertices   $u,v,w \in V$  such that $v$ is orthogonal to both $u$ and  $w$, then $u\sim w$, where $\sim$ is an equivalent relation on $G$  given by $u\sim w$ if     their open neighborhoods  coincide.   In \cite[Theorem 3.5]{DRJ}, these notions were  used, for the classical zero-divisor graph,  to characterize when the total quotient ring of a reduced ring $R$ is von Neumann regular.  Also,  \cite[Proposition 4.8]{DJF}  gives a similar result. It  was  shown that, when $\overline{\Gamma}(R)\neq \Gamma(R)$,  $\overline{\Gamma}(R)$ is complemented is a sufficient condition so that the total quotient ring  of  $R$ be zero-dimensional.  But, when checking the proof it seems that the proof holds true for $R$ with  $girth(\overline{\Gamma}(R))=4$.  In this paper, using a new treatment,  we prove that \cite[Proposition 4.8]{DJF} still holds true with any assumption (see Theorem \ref{0dim2}). So, in this paper, we continue the investigation begun in \cite{DJF} of the extended zero-divisor graph $\overline{\Gamma}(R)$ of a commutative ring $R$. Namely,  we study when   $\overline{\Gamma}(R)$ is complemented or uniquely complemented.\\

This article is organized as follows:\\

In Section 2,  we study when the extended zero-divisor graph of a commutative ring is complemented. We start by showing that,  if $\overline{\Gamma}(R)$ is complemented such that $|Z(R)|\geq 4$, then the ring $R$  has at most   one nonzero nilpotent element (see Theorem \ref{thmNil(R)} and Example \ref{exm_thm_Nil(R)}).  When $R$ is finite,  we have  the converse of Theorem \ref{thmNil(R)} (see Corollary \ref{conversenilcomp}). To show this, we give first a complete characterization for finite rings with complemented extended zero-divisor graphs (see Theorem \ref{them3}).\\

In Section 3, we show as a main result that,  for a  non reduced ring, if the extended zero-divisor graph is complemented,  then the complementeness and the uniquely complementeness properties coincide (see Theorem \ref{thmequivalentofuniqueness}). Also, if the extended zero-divisor graph  is complemented, then every orthogonal  of the nonzero nilpotent element is not an end (see Corollary \ref{vertexEnd}).\\

In Section 4, we show that the total quotient ring   $T(R)$ of $R$ is zero-dimensional once $\overline{\Gamma}(R)$ is complemented (see Theorem \ref{0dim2}). 
%Namely,  if $|Z(R)|=2$, then $R$ is isomorphic to $\Z_4$ or $\Z_2[X]/(X^2)$ and so $T(R)$ is zero-dimensional. If $|Z(R)|=3$, then $R$ is isomorphic to $\Z_2\times \Z_2$, $\Z_9$ or $\Z_3[X]/(X^2)$ and so, $\overline{\Gamma}(R)$ is complemented and $T(R)$ is zero-dimensional. Theorem \ref{0dim2} present the case when $|Z(R)|\geq 4$,  %. In fact, this result is already presented in \cite[Proposition 4.8]{DJF},  but the  proof is not correct. The converse of Theorem \ref{0dim2} is not always true, to see that we give   a counterexample in Example \ref{counterexample1}. Also,  Theorem \ref{0dim2}
%which leads us to show, 
This results  helps to show that  when 	$\overline{\Gamma}(R)$ is complemented,  every non nilpotent element has an orthogonal which is not nilpotent (see Theorem \ref{0dim1}). At the end of this section we prove that for any ring $R$ such that $|\nil(R)|=2$,  $R$ is not local or $\overline{\Gamma}(R)$ is not  complemented (see Proposition \ref{localorcomplemented}). \\

Finally, Section 5 is devoted to  study  when the extended zero-divisor graph of a direct  product of two rings  and the extended zero-divisor graph of an idealization of an $R-$module $B$  are complemented (see Theorems \ref{products1}, \ref{products2}, \ref{products3} and Proposition \ref{idealization}).

	%%%%%%%%%%%%%%%%%%%%%%%%%%%%%%%%%%%%%%%%%%%%%%%%%%%%%%%%%
	%%%%%%%%%%%%%%%%%%%%%%%%%%%%%%%%%%%%%%%%%%%%%%%%%%%%%%%%%
	%%%%%%%%%%%%%%%%%%%%%%%%%%%%%%%%%%%%%%%%%%%%%%%%%%%%%%%%%
	%%%Section 2
	%%%%%%%%%%%%%%%%%%%%%%%%%%%%%%%%%%%%%%%%%%%%%%%%%%%%%%%%%
	%%%%%%%%%%%%%%%%%%%%%%%%%%%%%%%%%%%%%%%%%%%%%%%%%%%%%%%%%
	%%%%%%%%%%%%%%%%%%%%%%%%%%%%%%%%%%%%%%%%%%%%%%%%%%%%%%%%%

	\section{When the extended  zero-divisor graph of  a commutative ring is complemented }
	In this section we study when the extended zero-divisor graph of a commutative ring is complemented. We start by showing that the ring $R$ will have at most  one nonzero nilpotent element  if $\overline{\Gamma}(R)$ is complemented. But first, we need the following lemmas which will be very useful through this paper.

	\begin{lem}\label{lem_indexofnilpotency}
		Let $R$ be a non reduced ring. If $\overline{\Gamma}(R)$ is complemented, then  every nonzero nilpotent element has index $2$.
	\end{lem}
	\pr
	Assume that $\nil(R)\neq \{0\}.$ Let $x\in \nil(R)$ such that $n_x\geq 3.$ Let $z\in Z(R)$ such that $z$ is adjacent to $x.$ If $x^{n_x-1}\neq z.$ Then,  $x^{n_x-1}$ is adjacent to both $z$ and $x.$ Otherwise, we can  easily see that $x^{n_x-1}+x$ is adjacent to both $x^{n_x-1}$ and $x$. Hence, $\overline{\Gamma}(R)$ is not complemented.\cqfd

\begin{figure}[ht]
		\centering
		\includegraphics[scale=0.5]{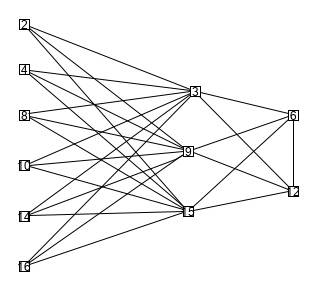}
		\caption{$\overline{\Gamma}(\Z_{18})$}
		\label{fig04}
	\end{figure}

	Notice that the converse of this lemma does not hold in general since, for instance, the extended zero-divisor    $\overline{\Gamma}(\Z_{18})$  illustrated in Figure \ref{fig04} is not complemented (since $ \overline{6}$ - for example - has no orthogonal element) even if the index of nilpotence  of every nilpotent element  is $2$. 
	
	\begin{exm}
		\begin{enumerate}
			\item Let $p$ be a prime number and $n$ be a positive  integer. Then, $\overline{\Gamma}(\Z_{p^n})$ is complemented  if and only if $n=2$ and $p=3$.  
			\item Consider   $\R[X,Y]/(X^3,XY^3)$.  The index of nilpotence of  $\overline{X}$ is $3$,  so $\overline{\Gamma}(\R[X,Y]/(X^3,XY^3))$ is not  complemented.  
		\end{enumerate}
	\end{exm}

	\begin{lem}\label{lemNil(R)}
		Let $R$ be a ring such that $|Z(R)^*|\geq 3$ and $\overline{\Gamma}(R)$ is complemented. Then, 
		\begin{enumerate}
			\item for every $\alpha\in \nil(R)^*$, $2\alpha= 0$, and 
			\item for every $\alpha\in \nil(R)^*$, if there exists $\beta\in Z(R)^*$ such that $\beta\perp\alpha$, then $\beta\notin \nil(R)$. 
		\end{enumerate} 
	\end{lem} 
	\pr 
	1. Assume that there exists $\alpha \in \nil(R)^*$ such that $2\alpha\neq 0$. Since $\overline{\Gamma}(R)$ is complemented, $\alpha$ is adjacent to $(-\alpha)$ on the other hand we have  $|Z(R)^*|\geq 3$, then every vertex $z\in Z(R)^*\setminus \{\alpha, -\alpha\}$ that is  adjacent to $\alpha$ is also adjacent  to $(-\alpha)$ and so $\alpha- z$ is a part of a triangle,  a contradiction.\\
	2. Let $\alpha\in \nil(R)^*$ such that $\alpha\perp \beta$ with $\beta\in Z(R)^*$. If $\beta\in \nil(R)^*$, then $\alpha+\beta\neq 0$, otherwise $\alpha=(-\beta)$ and so $2\alpha=0=\alpha-\beta$ which implies that $\alpha=\beta$,  a contradiction since $\alpha\perp \beta$. Thus, $\alpha+\beta$ is adjacent to both $\alpha$ and $\beta$ since $\alpha$ and $\beta$ are adjacent. And, by Lemma \ref{lem_indexofnilpotency},  $\beta^2=\alpha^2=0$. So,  $\alpha$ and  $\beta$ are not orthogonal,  a contradiction. Hence, $\beta\notin \nil(R)^*$.
	\cqfd
	
	%If  $R$ is reduced, then  in this case  $\Gamma(R)= \overline{\Gamma}(R)$. So, in the following results we treat the case when $R$ is not reduced. 
	Notice that, if  $|Z(R)|=2$, then $R$ is isomorphic to  $\Z_4$ or $\Z_2[X]/(X^2)$ and so $\nil(R)=\{0,a\}$ for some $0\neq a\in R$. In this case $\overline{\Gamma}(R)$ is not complemented. If $|Z(R)|=3$, then  $R$ is isomorphic to $\Z_9$ or $\Z_3[X]/(X^2)$ and so $\nil(R)=\{0,a,-a\}=Z(R)$  for some $0\neq a\in R$ and in this case $\overline{\Gamma}(R)$ is  complemented. When $|Z(R)|> 3$, we have the following result.

	\begin{thm}\label{thmNil(R)}
		Let $R$ be a ring such that $|Z(R)|\geq 4$.	If $\overline{\Gamma}(R)$ is  complemented, then  $|\nil(R)|\leq 2$.
	\end{thm}

	\pr
	Assume that there exist  $a, b\in \nil(R)^*$ such that $a\neq b$. Then, $a+b\in \nil(R)^*$, by Lemma \ref{lemNil(R)}.  Let $x,y,z\in Z(R)\setminus \nil(R)$ such that $x\perp a$, $y\perp b$ and $z\perp a+b$.  Let  $n$ be a positive integer such that $z^n(a+b)=0$,  we have the two following cases:\\
	\textbf{Case} $ab\neq 0$: We have $z^n(a+b)=0$,  then $z^nab=-z^nb^2=0$ by Lemma \ref{lem_indexofnilpotency},  and so   $ab$ is adjacent to both $z$ and $a+b$ ($ab\neq z$ since $ab \in \nil(R)^*$ and also $ab\neq a+b$). Then,  $z-(a+b)$ is a part of a triangle,  a contradiction.\\
	\textbf{Case} $ab=0$: If $z^na=0$, then $a$ is adjacent to both $z$ and $a+b$, a contradiction. Then, $z^na\neq0$. If $z^na\neq a$, then $z^na$ is adjacent to both $a$ and $x$ and so $x-a$ is a part of a triangle,   a contradiction. Otherwise,   we have  $z^n(a+b)=0$ and $b\in \nil(R)^*$, then  $z^na=-z^nb=z^nb$,  and so $z^na=a=z^nb$  is adjacent to both $b$ and $y$. Then,  $b-y$ is a part of a triangle, also a contradiction.
	\cqfd 

	\begin{exm}\label{exm_thm_Nil(R)}
		\begin{enumerate}
			\item Consider  the ring $R=D\times \Z_2[X]/(X^2)$,  where $D$ is an integral domain. Then,  $\nil(R)=\{(0,\bar{0}),(0,\bar{X})\}$   and its  extended zero-divisor graph is illustrated in Figure \ref{fig01}. Thus, $\overline{\Gamma}(R)$  is  a complete bipartite graph and hence it is  complemented.
			\begin{figure}[ht]
				\centering
				\includegraphics[scale=0.5]{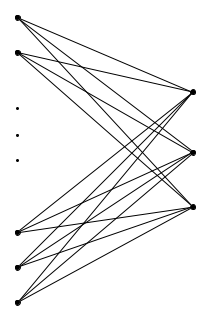}
				\caption{$\overline{\Gamma}(D\times \Z_2[X]/(X^2))$}
				\label{fig01}
			\end{figure}
			\item For the ring $R=\Z_2\times\Z_2\times\Z_4$, we have   $\nil(R)=\{(\bar{0},\bar{0},\bar{0}),(\bar{0},\bar{0},\bar{2})\}$. The extended zero-divisor graph of  this ring is illustrated in Figure \ref{fig02}.
			\begin{figure}[ht]
				\centering
				\includegraphics[scale=0.5]{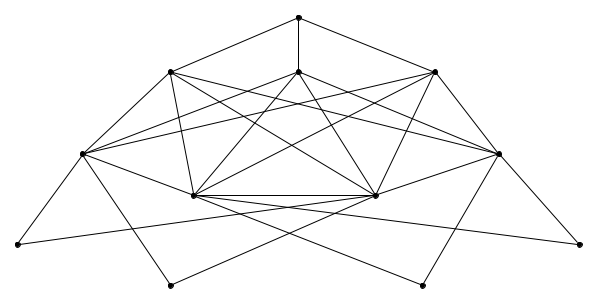}
				\caption{$\overline{\Gamma}(\Z_2\times\Z_2\times\Z_4)$}
				\label{fig02}
			\end{figure}
			Hence,  $\overline{\Gamma}(R)$ is complemented.
			
			\item For the ring $\Z_2[X,Y]/(X^3,XY)$, we have $\nil(\Z_2[X,Y]/(X^3,XY))=\{\bar{0}, \bar{X}, \bar{X}^2, \bar{X}+\bar{X}^2\}$. The extended zero-divisor of this ring is illustrated in Figure \ref{fig03}.
			\begin{figure}[ht]
				\centering
				\includegraphics[scale=0.5]{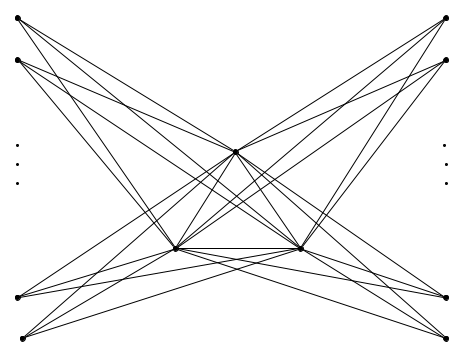}
				\caption{$\overline{\Gamma}(\Z_2[X,Y]/(X^3,XY))$}
				\label{fig03}
			\end{figure}

			Thus, $\overline{\Gamma}(\Z_2[X,Y]/(X^3,XY))$ is not complemented. Namely, $\bar{X}+\bar{Y}$ has no orthogonal element.
		\end{enumerate}
	\end{exm}
	When $R$ is finite,  we have  the converse of Theorem \ref{thmNil(R)} (see Corollary \ref{conversenilcomp}). To show this result, we give the following one.

	\begin{thm}\label{them3}
		Let $R$ be a finite ring such that $\Gamma(R)\neq \overline{\Gamma}(R)$. $\overline{\Gamma}(R)$ is complemented if and only if $R \cong B\times A_1\times\dots \times A_n $ where  $B\cong \Z_4$ or $\Z_2[X]/(X^2)$ and  $A_1,\dots , A_n$ are finite fields.
	\end{thm}

	\pr 
	$\Leftarrow$)%It follows from Lemma \ref{z4product1}.\\ 
	We can give a direct proof here. But, this can be proven easily, by induction, using Theorem  \ref{products1} and  \ref{products2} given in Section 5.\\  
	$\Rightarrow$) Since $R$ is a finite ring. Then,  by \cite[Theorem 87]{AM},  $R\cong A_1\times\cdots\times A_n$ such that $A_i$ is a finite local  ring for all $i\in\{1,\ldots,n\}$. Then,  for  all $i\in\{1,\ldots,n\}$, $Z(A_i)=\nil(A_i)$. By Theorem \ref{thmNil(R)}, $|\nil(R)|\leq 2$. Since $\overline{\Gamma}(R)\neq \Gamma(R)$, $|Nil(R)|=2$ and so one of the $A_i's$ is isomorphic to $\Z_4$ or $\Z_2[X]/(X^2)$  and the other rings are finite fields.
	\cqfd

	\begin{cor}\label{produitZ_4complemented}
		Let $n\in {\N}^* $ such that $\Gamma(\Z_n)\neq \overline{\Gamma}(\Z_n)$. Then,  $\overline{\Gamma}(\Z_n)$ is  complemented  if and only if  $n=2^2p_1\dots p_r$ with  $p_1,\dots , p_r$ are distinct prime  numbers and $r\geq 1$ is a positive integer.
	\end{cor}
	
Now, let us prove  the converse of Theorem \ref{thmNil(R)} in the case of a finite ring. 
	\begin{cor}\label{conversenilcomp}
		Let $R$ be a finite ring such that $\Gamma(R)\neq \overline{\Gamma}(R)$. If $Nil(R)=\{0,a\}$ for some $a\in R^*$, then $\overline{\Gamma}(R)$ is complemented.\end{cor}
	\pr
	Since $R$ is a finite ring,  by \cite[Theorem 87]{AM},  $R\cong A_1\times\cdots\times A_n$ such that $A_i$ is a finite local  ring for all $i\in\{1,\ldots,n\}$. Then,  for  all $i\in\{1,\ldots,n\}$, $Z(A_i)=\nil(A_i)$. Then, $R$ is either indecomposable ring such that $R$ is a finite local ring or $R\cong A_1\times\cdots\times A_n$ such that $A_i$ is a finite local  ring for all $i\in\{1,\ldots,n\}$, where $n\geq 2$. If $R$ is indecomposable, then, using the fact that $|Nil(R)|=2=|Z(R)|$, $R\cong \Z_4$. Then, this contradicts the fact that $\Gamma(R)\neq \overline{\Gamma}(R)$. Thus,  $R\cong A_1\times\cdots\times A_n$ such that $Z(A_i)=\nil(A_i)$ for every $i\in\{1,\ldots,n\}$ and $n\geq 2$. Since $|Nil(R)|=2$,  one of the $A_i's$  is isomorphic to $\Z_4$ or $\Z_2[X]/(X^2)$  and the other rings are integral domains. Then, by Theorem \ref{them3}, $\overline{\Gamma}(R)$ is complemented. 
	\cqfd

The authors are not able to prove the equivalence of Theorem \ref{them3} for infinite ring. We let it then as an open interesting  question.  
	
%In the previous result, we proved that if the ring $R$ is isomorphic to $B \times A_1\times\dots\times A_n$ such that $B\cong \Z_4$ or $ \Z_2[X]/(X^2)$ and  $A_1,\ldots,A_n$ are integral domains, then $\overline{\Gamma}(R)$ is complemented. If the ring is finite, then we have the equivalence. Thus,  there is a natural question that arises in this context which is: Have we this equivalence also for the infinite rings or there exist some infinite rings,  which are not isomorphic to $B \times A_1\times\dots\times A_n$, such that its associate extended zero-divisor graph is complemented. If the answer is positive, what is the form of these kind of rings?

	\section{Complementeness and uniquely complementeness properties coincide in the extended zero-divisor graphs}
	
	In \cite[Theorem  3.5]{DRJ}, it was shown that, when $R$ is reduced,  $\Gamma(R)(=\overline{\Gamma}(R))$ is uniquely complemented if and only if $\Gamma(R)$ is complemented if and only if $T(R)$ is von Neumann regular. The main result of this section generalizes \cite[Theorem  3.5]{DRJ}. Namely, shows that, when $R$ is not reduced, the complementeness and the uniquely complementeness properties coincide. To show this, we first prove the following lemma. 
	
	\begin{lem}\label{propertyofothogonality}
		Let $R$ be a ring and  $a, b, c \in Z(R)\setminus \nil(R)$ such that $a\perp b$ and $a\perp c$, then $b\sim c$.
	\end{lem}
	\pr
	We have $a^{n_1} b^{m_1}=a^{n_2} c^{m_2}=0$  for some $n_1, m_1, n_2,m_2 \in \N$. We first show that $b$ and $c$ are not adjacent, that is $b^\alpha c^\beta\neq 0$ for every $\alpha,\beta\in \N$. If $b^\alpha c^\beta =0$ for some $\alpha, \beta\in \N$. Then, $b=c$ or $a=c$ (since $a\perp b$ and $a\perp c$) and so $b\in \nil(R)$ or $a\in \nil(R)$, a contradiction. Thus, $b$ and $c$ are not adjacent. Now, let prove that $N(b)=N(c)$. Let $d\in N(b)$, then $d^nb^m=0$ with $d^n\neq 0$ for some $n,m\in \N$. Thus, $(d^nc^{m_2})a^{n_2}=d^n(c^{m_2}a^{n_2})=0$ and $(d^nc^{m_2})b^{m}=(d^nb^{m})c^{m_2}=0$. So, $d^nc^{m_2}=0$, otherwise $d^nc^{m_2}$ is adjacent to both $a$ and $b$ (and $d^nc^{m_2}\neq a$ or $b$ since $a,b\notin \nil(R)$) which contradict the fact that $a\perp b$. Hence, $d\in N(c)$. Similarly, we show the other inclusion. Thus, $b\sim c$.
	\cqfd

	Now, we are ready to prove the main result of this section.

	\begin{thm}\label{thmequivalentofuniqueness}
		Let $R$ be a ring such that $\overline{\Gamma}(R)\neq\Gamma(R)$. Then, $\overline{\Gamma}(R)$ is uniquely complemented if and only if $\overline{\Gamma}(R)$ is complemented. 
	\end{thm}
	\pr
	$\Rightarrow)$ By definition of  uniquely complemented.\\
	$\Leftarrow)$ Suppose that $\overline{\Gamma}(R)$ is  complemented, then, by Theorem 2.5, $\nil(R)=\{0,\alpha\}$ for some $0\neq \alpha \in R$ and so by   Lemma \ref{propertyofothogonality}, we have just to prove that for every $b,c\in Z(R)^*$, if  $\alpha\perp b$ and $\alpha\perp c$, then $b \sim  c$ and if  $\alpha\perp c$ and $b\perp c$, then $\alpha \sim b$. Let prove the first one. So, suppose by contradiction  that  there exist $b, c\in Z(R)^*$ such that $\alpha \perp b$ and $\alpha \perp c$ but $b\nsim c$. Then, there exists $x\in N(c)\setminus N(b)$ and so $x^{n_1}c^{m_1}=0$ for some $n_1,m_1\in \N^*$ and $x^nb^m\neq 0$ for every $n,m\in \N^*$. Assume that $xb\neq c$, then $(xb)^{n_1} c^{m_1}=0$ and so $xb$ and $c$ are adjacent. On the other hand $\alpha$ and $b$ are adjacent which implies that $\alpha b^t=0$ for some $t\in \N$, then $\alpha (xb)^t =0$ and so $xb$ is adjacent to both  $\alpha$ and  $c$ a contradiction since $\alpha\perp c$. Then, $xb=c$ and so  $x^{n_1}c^{m_1}=0$ implies $x^{n_1}(xb)^{m_1}=x^{n_1+m_1} b^{m_1}=0$ a contradiction. Now, we prove the second one, that is,  for every $b,c\in Z(R)^*$ such that $\alpha \perp c$ and $b\perp c$, then $\alpha \sim b$. Assume that $\alpha \perp c$ and $b\perp c$, then $\alpha c^{m_1}=b^{n_1} c^{m_2}=0$ for some $n_1,m_1,m_2\in \N^*$. If $\alpha$ is adjacent to $b$, then  $b$ is  adjacent to both $\alpha$ and $c$, a contradiction since $\alpha\perp c$. Thus, $\alpha$ is not adjacent to $b$. Let $d\in N(\alpha)$, then $d^n\alpha =0$ for some $n\in \N^*$ and so $(d^nb^{n_1})c^{m_2}=d^n(b^{n_1}c^{m_2})=0$ and $(d^nb^{n_1})\alpha=b^{n_1}d^n\alpha=0$. If $d^nb^{n_1}\in Z(R)\setminus \nil(R)$, $d^nb^{n_1}\neq \alpha$ and $d^nb^{n_1}\neq c$ and so $d^nb^{n_1}$ is adjacent to both $c$ and $\alpha$, a contradiction (since $\alpha\perp c$). Thus, $d^nb^{n_1}\in \nil(R)$. We  have then two cases to discuss:\\
	\textbf{Case 1}. $d^nb^{n_1}=0$. Then, $d$ is adjacent to $b$ since $d\neq b$ and $d,b\in Z(R)\setminus \nil(R)$. Thus, $d\in N(b)$.\\
	\textbf{Case 2}. $d^nb^{n_1}=\alpha$. Then, $d^{2n}b^{2n_1}=\alpha^2=0$ with $d^{2n}\neq 0$, $b^{2n_1}\neq 0$ and $d\neq b$ (since $b,d\in Z(R)\setminus \nil(R)$). Hence, $d\in N(b)$.\\  Then, $\alpha\sim b$ and  hence  $\overline{\Gamma}(R)$ is uniquely complemented.
	\cqfd

	\begin{cor}
		Let $R$ be a ring such that $\Gamma(R)\neq \overline{\Gamma}(R)$ and $\overline{\Gamma}(R)$ is complemented. Then, for every orthogonal $b\in Z(R)^*$ of the nonzero nilpotent element $\alpha$, we have $b\sim \alpha+b$. 
	\end{cor}
	\pr
	Assume that $\Gamma(R)\neq \overline{\Gamma}(R)$ and   $\overline{\Gamma}(R)$ is complemented, then by Theorem \ref{lemNil(R)}, $Nil(R)=\{0,\alpha\}$ for some $0\neq \alpha\in R$.
	Let $b\in Z(R)^*-\{\alpha\}$ such that $\alpha\perp b$, so $\alpha b^n=0$ for some positive integer $n$ and there is no vertex adjacent to both $\alpha$ and $b$. First, let prove that $\alpha\perp (\alpha+b)$. We have $\alpha(\alpha+b)^n=\alpha(b^n+n\alpha b^{n-1}+\ldots +\alpha^n)=\alpha b^n=0$ and  since $\alpha+b\neq \alpha$ and $(\alpha+b)^n\neq 0$ (because,  $b\notin \nil(R)$), then $\alpha$ and $\alpha+b$ are adjacent. Now, assume that there exists $c$ that is adjacent to both $\alpha$ and $\alpha+b$, then $c^{n_1}\alpha=0=c^{n_1}(\alpha+b)^{m_1}=m_1c^{n_1}\alpha b^{m_1-1} +c^{n_1}b^{m_1}=0+c^{n_1}b^{m_1}$ and so $c$ is adjacent to $b$,  a contradiction since  $\alpha\perp b$. Hence, $\alpha\perp (\alpha+b)$. Thus, by Theorem  \ref{thmequivalentofuniqueness},  $b\sim \alpha+b$.
	\cqfd

	\section{Complemented extended zero-divisor graphs and zero-dimensional rings}
	
	If $|Z(R)|=2$, then $R$ is isomorphic to $\Z_4$ or $\Z_2[X]/(X^2)$ and so $T(R)$ is zero-dimensional. If $|Z(R)|=3$, then $R$ is isomorphic to $\Z_2\times \Z_2$, $\Z_9$ or $\Z_3[X]/(X^2)$ and so, $\overline{\Gamma}(R)$ is complemented and $T(R)$ is zero-dimensional. In this section, we show that $T(R)$ is zero-dimensional once $\overline{\Gamma}(R)$ is complemented. In fact, this result was already given in \cite[Proposition 4.8]{DJF}.  But, in the third line of the proof, \cite[Corollary 3.4]{DJF} is used to show that an element $z_0$ is not nilpotent. This means that like we have supposed that the girth of $\overline{\Gamma}(R)$ is not $3$. But, we can give examples of $\overline{\Gamma}(R)$ which are complemented  with girth equal to $3$. For this consider $\overline{\Gamma}(\Z_2\times \Z_2 \times \Z_4)$ (see Figure \ref{fig02}). Now, using a new way, we show that \cite[Proposition 4.8]{DJF} holds true. To show that, we need the following lemma.

	\begin{lem}\label{uniquelycomplemented}
		Let $R$ be a ring such that $\Gamma(R) \neq \overline{\Gamma}(R)$. If $\Gamma(R)$ is uniquely complemented, then $\overline{\Gamma}(R)$ is not complemented.
	\end{lem}
	\pr
	If $\Gamma(R)$ is uniquely complemented, then, using \cite[Theorem 3.9]{DRJ}, $\Gamma(R)$ is a star graph and so $\overline{\Gamma}(R)$ is not complemented.
	\cqfd
	
	Using the previous lemma, we get the main result of this section. 
	
	\begin{thm}\label{0dim2}
		Let $R$ be a ring such that $|Z(R)^*|\geq 3$ and $\Gamma(R)\neq \overline{\Gamma}(R)$. If $\overline{\Gamma}(R)$ is complemented then $T(R)$ is zero-dimensional.
	\end{thm}
	
	\pr
	There are two cases to discuss:\\
	\textbf{Case 1.} For every $x\in Z(R)^*$,  $x^{\perp}\cap (Z(R)\setminus \nil(R))\neq \emptyset$. In this case, we show that for every $ \frac{x_1}{x_2}$  in $T(R)$, there exists $\frac{m_1}{m_2}\in T(R)$ such that $\frac{x_1}{x_2}+\frac{m_1}{m_2}$ is unit and $\frac{x_1}{x_2}\frac{m_1}{m_2}$ is nilpotent. Then, let $\frac{x_1}{x_2}$  in $T(R)$. So,  there are three sub-cases to discuss:\\
	\textbf{Sub-case 1.} Assume that $x_1 \in R\setminus Z(R)$. Since $\overline{\Gamma}(R)$ is complemented and $|Z(R)^*|\geq 3$, $|\nil(R)^*|=1$. We denote by $\alpha$ the nilpotent element of $R$.  Using Lemma \ref{lemNil(R)}, we have $\alpha^2=2\alpha=0$. It is clear that $\frac{x_1}{x_2}\frac{\alpha}{x_2}$ is nilpotent, also $\frac{x_1}{x_2}+\frac{\alpha}{x_2}$ is unit since $(x_1+\alpha)^2={x_1}^2\notin Z(R)$.\\
	\textbf{Sub-case 2.} Assume that $x_1 =\alpha$. We have $\frac{x_1}{x_2}\frac{1}{x_2}$ is nilpotent,  also  $\frac{x_1}{x_2}+\frac{1}{x_2}$ is unit since $(x_1+1)^2=1\notin Z(R)$. \\
	\textbf{Sub-case 3.} Assume that $x_1 \in Z(R)\setminus \nil(R)$. Then,  there exists $m_1\in {x_1}^{\perp}\cap (Z(R)\setminus \nil(R))$. We have $\frac{x_1}{x_2}\frac{m_1}{x_2}$ is nilpotent (since $x_1$ and $m_1$ are adjacent), it remains to show that $\frac{x_1}{x_2}+\frac{m_1}{x_2}$ is unit, for that it is sufficient to prove that $x_1+m_1$ does not belong to $Z(R)$.	Let $z\in R$ such that $z(x_1+m_1)=0$.   We have $x_1m_1$ is nilpotent (since $x_1$ and $m_1$ are adjacent), then there are  the two following  sub-subcases to discuss.\\
	\textbf{Sub-subcase 1.} If $x_1m_1=0$. We have $z(x_1+m_1)=0$, then $zx_1m_1+z{m_1}^2=0$ and $z{x_1}^2+zx_1m_1=0$, so  $z{m_1}^2=0$ and $z{x_1}^2=0$. Then, $z\neq x_1$ and $z\neq m_1$ since $x_1$ and $m_1$ are not nilpotent, then $z$ is adjacent to both $x_1$ and $m_1$, a contradiction (since  $x_1$ and $m_1$ are orthogonal). Then $z=0$.\\
	\textbf{Sub-subcase 2.} If $x_1m_1=\alpha$. We have $z{x_1}^2+zx_1m_1=0$ and $zx_1m_1+z{m_1}^2=0$, then  $z{x_1}^2+z\alpha=0$ and $z\alpha+z{m_1}^2=0$. 
	$z\alpha=0$, otherwise  $z{x_1}^2+z\alpha=0$ and $z\alpha+z{m_1}^2=0$ imply that  $z\alpha{x_1}^2=0$ and $z\alpha{m_1}^2=0$. Then, $z\alpha\neq x_1$, $z\alpha\neq m_1$ since $x_1$ and $m_1$ are not nilpotent and so    $z\alpha$ is adjacent to both $x_1$ and $m_1$, a contradiction (since  $x_1$ and $m_1$ are orthogonal). Then $z\alpha=0$ which implies that 
	$z{x_1}^2=0=z{m_1}^2=0$. Then, $z\neq x_1$ and $z\neq x_1$ since $x_1$ and $m_1$ are not nilpotent, then $z$ is adjacent to both $x_1$ and $m_1$, a contradiction (since  $x_1$ and $m_1$ are orthogonal). Then,  $z=0$.\\  Hence, $\frac{x_1}{x_2}+\frac{m_1}{x_2}$ is unit.\\ Then, $T(R)$ is $\pi$-regular and so it is zero-dimensional.\\
	\textbf{Case 2.} There exists $x\in Z(R)^*$,  $x^{\perp}=\{\alpha\}\subset \nil(R)=\{0,\alpha\}$. In this case, one can show that $N_{\overline{\Gamma}(R)}(x)=\{\alpha\}$. Otherwise, there exist $s$ and $t$ in $N_{\overline{\Gamma}(R)}(x)$ such that $s$ and $t$ are adjacent. Since $\{\alpha\}=x^{\perp}$, $s\alpha=\alpha$. Also, since $s$ and $t$ are adjacent ,$st=\alpha$ or $st=0$. Then, $s\alpha t=0$,  which implies that $\alpha t=0$, then $t$ and $\alpha$ are adjacent, which is absurd with the fact that $\{\alpha\}=x^{\perp}$. Thus, $N_{\overline{\Gamma}(R)}(x)=\{\alpha\}$. On the other hand, $\Gamma(R)$ is not uniquely complemented (using Lemma \ref{uniquelycomplemented}). Then, in this case, there are two cases to discuss:\\
	\textbf{Sub-case 1.} $\Gamma(R)$ is complemented. Since $\Gamma(R)$ is not uniquely complemented, by \cite[Theorem 3.14]{DRJ}, $R\simeq \Z_4\times D$ or $\Z_2[X]/(X^2)\times D$ such that $D$ is an integral domain. Then, $T(R)$ is zero-dimensional.\\
	\textbf{Sub-case 2.} $\Gamma(R)$ is not complemented. Then, there exists $b\in Z(R)^*$ such that $b$ has an orthogonal in $\overline{\Gamma}(R)$ and not in $\Gamma(R)$ (one can see that $b\neq \alpha$ since $\alpha$ has $x$ as an orthogonal in $\Gamma(R)$). Then, there exists $t\in b^{\perp}$ such that $bt=\alpha$. One can show that $t\neq \alpha$. Otherwise, $b\alpha=\alpha$ and so for every $n\in \N^*$,  $b^n\alpha=\alpha\neq 0$ and so $b$ and $\alpha$ are not adjacent in $\overline{\Gamma}(R)$, which is absurd with the fact that $t=\alpha$ is orthogonal to $b$. There are two cases to discuss:\\
	\textbf{Sub-subcase 1.}
	$b\alpha\neq 0$. Then, $b\alpha=\alpha$. Since for every $z\in N_{\overline{\Gamma}(R)}(b)$, $zb=\alpha$ or $zb=0$, $z(b\alpha)=0$ for every  $z\in N_{\overline{\Gamma}(R)}(b)$, then $z\alpha=0$ for every $z\in N_{\overline{\Gamma}(R)}(b)$ 
	(since $b\alpha=\alpha$). In this case, we show that ${(bx)}^{\perp}=\emptyset$, therefore, we determine firstly $N_{\overline{\Gamma}(R)}(bx)$. Let $h\in N_{\overline{\Gamma}(R)}(bx)$, then there exists $n,m$ in $\N^*$ such that $(bx)^nh^m=0$ with $(bx)^n\neq 0$ and  $h^m\neq 0$. Then, $(bh)^nx^n=0$ (resp., $(bh)^mx^n=0$), if $n\geq m$ (resp., $m\geq n$). If $(bh)^n\neq 0$, then $bh=\alpha$ since $N_{\overline{\Gamma}(R)}(x)=\{\alpha\}$. Then, $h=\alpha$ or $h\in N_{\overline{\Gamma}(R)}(b)$. If $(bh)^n=0$, then $h=\alpha$ or $h\in N_{\overline{\Gamma}(R)}(b)$. Thus, $N_{\overline{\Gamma}(R)}(bx)=N_{\overline{\Gamma}(R)}(b)\cup \{\alpha\}$. Thus, ${(bx)}^{\perp}=\emptyset$ since $z\alpha=0$ for every $z\in N_{\overline{\Gamma}(R)}(b)$, which is absurd with the fact that $\overline{\Gamma}(R)$ is complemented.\\
	\textbf{Sub-subcase 2.} $b\alpha=0$. Then,  $t\alpha \neq 0$ since $t\in b^{\perp}$, which implies that $t\alpha=\alpha$. As in the previous case, we show that for every $z\in N_{\overline{\Gamma}(R)}(t)$, $z\alpha=0$. In this case, we show that ${(tx)}^{\perp}=\emptyset$, therefore, we determine firstly $N_{\overline{\Gamma}(R)}(tx)$. 
	Let $h$ in $N_{\overline{\Gamma}(R)}(tx)$,  then there exists $n,m$ in $\N^*$ such that $(tx)^nh^m=0$ with $(tx)^n\neq 0$ and $h^m\neq 0$. Then, $(th)^nx^n=0$ (resp. $(th)^mx^n=0$), if $n\geq m$ (resp., $m\geq n$). If $(th)^n\neq 0$, then $th=\alpha$ since $N_{\overline{\Gamma}(R)}(x)=\{\alpha\}$. Then, $h=\alpha$ or $h\in N_{\overline{\Gamma}(R)}(t)$. If $(th)^n=0$, then $h=\alpha$ or $h\in N_{\overline{\Gamma}(R)}(t)$. Thus, $N_{\overline{\Gamma}(R)}(tx)=N_{\overline{\Gamma}(R)}(t)\cup \{\alpha\}$. Thus, ${(tx)}^{\perp}=\emptyset$ since for every $z\in N_{\overline{\Gamma}(R)}(t)$, $z\alpha=0$, which is absurd with the fact that $\overline{\Gamma}(R)$ is complemented.
	\cqfd

	$T(R)$ is zero dimensional does not imply necessarily that $\overline{\Gamma}(R)$ is complemented, as example, we have the following one.  
	
	\begin{exm}\label{counterexample1}
		$T(\Z_{16})$ is zero dimensional and $\overline{\Gamma}(\Z_{16} )$ is not complemented.
	\end{exm}
	
	\pr
	Since $\Z_{16}$  is zero-dimensional, $T(\Z_{16})\simeq \Z_{16}$ is zero-dimensional. On the other hand, $\overline{\Gamma}(\Z_{16})$ is not complemented, indeed:  $\overline{2}$ is a nilpotent element in $ \Z_{16}$ of index of nilpotence $4$. Then, by Lemma \ref{lem_indexofnilpotency},  $\overline{\Gamma}(\Z_{16})$ is not complemented. 
	\cqfd

	Theorem \ref{0dim2} leads us to show, when 	$\overline{\Gamma}(R)$ is complemented, that every non nilpotent element has an orthogonal which is not nilpotent, namely we have the following result.
	
	\begin{thm}\label{0dim1}
		Let $R$ be a ring such that $\overline{\Gamma}(R)$ is complemented and $|Z(R)^*|\geq 3$. Then, for all $x\in Z(R)\setminus \nil(R)$, $x^{\perp}\cap (Z(R)\setminus \nil(R))\neq \emptyset$.
	\end{thm}
	\pr
	Since $\overline{\Gamma}(R)$ is complemented and $|Z(R)^*|\geq 3$, $|\nil(R)^*|=1$ and $n_x\leq 2$ for every nilpotent element $x$. We denote by $\alpha$ the non zero nilpotent element of $R$.
	We suppose that there exists $x_1\in Z(R)\setminus \nil(R)$ such that
	${x_1}^{\perp}=\alpha$. On the other hand, by Theorem \ref{0dim2}, $T(R)$ is zero-dimensional, then, using \cite[Theorem 3.1, Theorem 3.2]{J}, there  exists $m_1$ in $R$ such that $\frac{x_1}{x_2}\frac{m_1}{m_2}$ is nilpotent and $\frac{x_1}{x_2}+\frac{m_1}{m_2}$ is unit for some $m_2,x_2$ in $R\setminus Z(R)$. $\frac{x_1}{x_2}\frac{m_1}{m_2}$ is nilpotent implies that $x_1m_1\nil(R)$. Then, there are two cases to discuss:\\
	\textbf{Case 1.} We assume that $x_1m_1=0$. Then, $x_1$ and $m_1$ are adjacent, if  $m_1$ is nilpotent, then ${m_1}^2=0$, so $m_1(x_1m_2+m_1x_2)=0$ which implies that $\frac{x_1}{x_2}+\frac{m_1}{m_2}$ is not unit, a contradiction.
	Otherwise, $m_1$ and $x_1$ are not orthogonal, then there exists $z\in Z(R)^*$ that adjacent to both $x_1$ and $m_1$. Then, $z^2{x_1}^2=z^2{m_1}^2=0$ (Since $zx_1$ and $zm_1$ are nilpotent). If $z^2{x_1}=0$ and $z^2{m_1}=0$, then $z^2({x_1}{m_2}+ {x_2}{m_1})=0$, then $\frac{x_1}{x_2}+\frac{m_1}{m_2}$ is not unit, a contradiction. Otherwise, $z^2{x_1}({x_1}{m_2}+ {x_2}{m_1})=0$ with $z^2{x_1}\neq 0$  or $z^2{m_1}({x_1}{m_2}+ {x_2}{m_1})=0$ with $z^2{x_1}\neq 0$, then  $\frac{x_1}{x_2}+\frac{m_1}{m_2}$ is not unit, a contradiction.\\
	\textbf{Case 2.} We assume that $x_1m_1=\alpha$. If $m_1$ is nilpotent, then $x_1$ and $m_1$ are adjacent, then there exists $n\in {\N}^*$ such that ${x_1}^nm_1=0$. we choose $\beta$ such that ${x_1}^\beta m_1=0$ and  ${x_1}^{(\beta-1)} m_1\neq 0$. We have ${x_1}^{(\beta-1)} m_1(x_1m_2+m_1x_2)$  then $\frac{x_1}{x_2}+\frac{m_1}{m_2}$ is not unit, a contradiction. If $m_1$ is not nilpotent, we have ${x_1}^2{m_1}^2=0$ with ${x_1}^2\neq 0$ and ${m_1}^2\neq 0$. Then $x_1$ and $m_1$ are adjacent and not orthogonal (since $m_1$ is not nilpotent). Then, there exists $z\in Z(R)^*$ that adjacent to both $x_1$ and $m_1$. Then, $z^2(x_1)^2=0$ and $z^2m_1=0$. If $z^2x_1=0$ and $z^2m_1=0$, then $z^2(x_1m_2+m_1x_2)=0$, then $\frac{x_1}{x_2}+\frac{m_1}{m_2}$ is not unit, a contradiction. Otherwise, $z^2x_1(x_1m_2+m_1x_2)=0$ if $z^2x_1\neq 0$ or $ z^2m_1(x_1m_2+m_1x_2)=0$ if $z^2m_1\neq 0$ and $z^2x_1=0$. then $\frac{x_1}{x_2}+\frac{m_1}{m_2}$ is not unit, a contradiction. Hence, for all $x\in Z(R)\setminus \nil(R)$, $x^{\perp}\cap (Z(R)\setminus \nil(R))\neq \emptyset.$ 
	\cqfd
	
	It was proven in \cite[Lemma 3.7]{DRJ} that if $\Gamma(R)$ is uniquely complemented and $|R|>9$, then any orthogonal of the nonzero nilpotent element of $R$ is an end, which is not the case for $\overline{\Gamma}(R)$. Namely, we have the following corollary:

	\begin{cor}\label{vertexEnd}
		Let $R$ be a ring  such that $\overline{\Gamma}(R)\neq \Gamma(R)$.  If $\overline{\Gamma}(R)$ is complemented, then every orthogonal of the nonzero nilpotent element is not an end.   
	\end{cor}

	Theorem \ref{0dim1} is useful to show that, when $Nil(R)=\{0,\alpha\}$ then either $R$ is not local or $\overline{\Gamma}(R)$ is not  complemented. Using Theorem \ref{conversenilcomp},  when $R$ is finite then 	$\overline{\Gamma}(R)$ is complemented. Then, in the following we can assume that $R$ is an infinite ring, namely we have the following result.
	
	\begin{pro}\label{localorcomplemented}
		Let $R$ be an infinite ring such that $Nil(R)=\{0,\alpha\}$. Then, either $R$ is not local or $\overline{\Gamma}(R)$ is not  complemented.	 
	\end{pro}
	\pr Assume that  $Nil(R)=\{0,\alpha\}$ and suppose that $R$ is local with maximal ideal $\ann(\alpha)$ and that $\overline{\Gamma}(R)$ is complemented. Let $x\in Z(R)\setminus\{0,\alpha\}$, then , by Theorem \ref{0dim1}, there exists $y\in Z(R)\setminus \{0,\alpha\}$ such that $x\perp y$. But, since $\ann(\alpha)$ is the maximal ideal of $R$, $x,y\in \ann(\alpha)$. So, $x-y$ is a part of a triangle a contradiction.
	\cqfd

	\section{When $\overline{\Gamma}(R_1\times R_2)$ and $\overline{\Gamma}(R(+)B)$ are complemented}

	In the following first results of this section, we investigate when the extended zero-divisor graph of the  product of two rings,  $R_1\times R_2$,  is complemented. Namely, we treat the three existed cases following the cardinality of $Z(R_2)$; $|Z(R_2)|=1$, $|Z(R_2)|=2$ and $|Z(R_2)|\geq 3$.\\ 
	For the case when $R_2$ is integral domain, we have the following theorem.

	\begin{thm}\label{products1}
		Let $R_1$ and $R_2$ be two rings such that $R_2$ is an integral domain. Then, $\overline{\Gamma}(R_1\times R_2)$ is complemented if and only if either $|Z(R_1)|=2$ or  $(\ \overline{\Gamma}(R_1)$ is complemented and $|\nil(R_1)|\leq 2 )$.
	\end{thm}
\pr 
$\Rightarrow)$ Assume that $\overline{\Gamma}(R_1\times R_2)$ is complemented and $|Z(R_1)|\neq 2$.  If $|\nil(R_1)|\geq 3$, then $|\nil(R_1\times R_2)|\geq 3$,  a contradiction, by Theorem \ref{thmNil(R)}. Now, suppose that $\overline{\Gamma}(R_1)$ is not complemented, then there exists $z\in Z(R_1)^*$ such that for every $x\in Z(R_1)^*$, $x$ is not an orthogonal of $z$. We have $(z,0)\in Z(R_1\times R_2)$. Let $(a,b)\in Z(R_1\times R_2)$ such that $(a,b)$ is adjacent to $(z,0)$. So, $(a,b)^n(z,0)^m=(0,0)$ for some $n,m\in \N^*$ with $(a,b)^n\neq (0,0)$ and  $(z,0)^m\neq (0,0)$, then $a^nz^m=0$ and so we  have three cases to discuss:\\ 
\textbf{Case 1.} If $a^n=0$ and $b\neq 0$, then  for every vertex $y$ adjacent to $z$, $(y,0)$ is adjacent to both $(a,b)$ and $(z,0)$.\\
\textbf{Case 2.} If $a^n\neq 0$ and $b\neq 0$, then $a$ is adjacent to $z$ and so there exists $x$ adjacent to both $a$ and $z$ since $z$ is not orthogonal to $a$. Thus, $(x,0)$ is adjacent to both $(z,0)$ and $(a,b)$.\\
\textbf{Case 3.} If $a^n\neq 0$ and $b=0$, then $(0,1)$ is adjacent to both $(a,b)$ and $(z,0)$.\\ 
Hence,  $(z,0)$ has no orthogonal in $\overline{\Gamma}(R_1\times R_2)$,  a contradiction.\\
$\Leftarrow)$ If $|\nil(R_1)|=1$, then $Z(R_1\times R_2)=(R_1\setminus Z(R_1)\times \{0\})\cup (Z(R_1)\times \{0\})\cup (Z(R_1)\times R_2^*)$. If $(a,b)\in R_1\setminus Z(R_1)\times \{0\}$, then $(a,0)\perp (0,1)$. If $(a,b)\in Z(R_1)\times \{0\}$, then $b=0$ and  $(a,0)\perp (c,1)$ with $c\in a^\perp$. If $(a,b)\in Z(R_1)\times R_2^*$, then $(a,b)\perp (c,0)$ with $c\in a^\perp$.\\ 
If $|\nil(R_1)|=2$, then $\nil(R_1)=\{0,\alpha\}$ for some $0\neq \alpha \in R_1$ and $Z(R_1\times R_2)=R_1\setminus Z(R_1)\times \{0\}\cup Z(R_1)\setminus \nil(R_1)\times R_2\cup \nil(R_1)\times R_2$. Let $(a,b)\in Z(R_1\times R_2)^*$. If $(a,b)\in R_1\setminus Z(R_1)\times \{0\}$, then $(a,0)\perp (0,1)$. If $(a,b)\in Z(R_1)\setminus \nil(R_1)\times R_2$, then if $b=0$, $(a,0)\perp (c,b)$ with $c\in a^\perp$, otherwise $(a,b)\perp (c,0)$ with $c\in a^\perp$. If $(a,b)\in \nil(R_1)\times R_2$, then if $(a,b)=(\alpha,0)$, $(\alpha,0)\perp (c,b')$ with $c\in \alpha^\perp$ and $b'\in R_2^*$, if $b\neq 0$ and $a=\alpha$, $(\alpha,b)\perp (c,0)$ with $c\in R_1\setminus Z(R_1)$, if $a=0$ and $b\neq 0$, $(0,b)\perp (c,0)$ with $c\in R_1\setminus Z(R_1)$.
\cqfd

Now, for the case when $|Z(R_2)|=2$, we have the following result.  

\begin{thm}\label{products2}
 Let $R_1$ and $R_2$ be two rings such that $|Z(R_2)|=2$. Then, $\overline{\Gamma}(R_1\times R_2)$ is complemented if and only if $\overline{\Gamma}(R_1)$ is complemented and $R_1$ is reduced.
\end{thm}
\pr
    $\Rightarrow)$ Assume that $\overline{\Gamma}(R_1\times R_2)$ is complemented.  We have $|\nil(R_2)|=|Z(R_2)|=2$, then if $R_1$ is not reduced, $|\nil(R_1\times R_2)|\geq 3$,  a contradiction (by Theorem \ref{thmNil(R)}). Now, suppose that $\overline{\Gamma}(R_1)$ is not complemented, then there exists $z\in Z(R_1)^*$ which has no orthogonal. We have $(z,0)\in  Z(R_1\times R_2)$. Suppose that there exists $(a,b)\in Z(R_1\times R_2)$ such that $(z,0)\perp (a,b)$. Then, $(a,b)^n(z,0)^m=(0,0)$ for some $n,m\in \N^*$ and so $a^nz^m=0$. Thus, we have two cases to discuss:\\
    \textbf{Case 1.}  If $a^n=0$, then $b^n\neq 0$ and so for every $y\in Z(R_1)^*$ such that $y$ is adjacent to $z$, $(y,0)$ is adjacent to both $(a,b)$ and $(z,0)$,  a contradiction.\\
     \textbf{Case 2.} If $a^n\neq 0$, then $a$ is adjacent to $z$ and so there exists $x\in Z(R_1)^*$ such that $x$ is adjacent to both $z$ and $a$. Thus, $(x,0)$ is adjacent to both $(a,b)$ and $(z,0)$,  a contradiction.\\
     Hence, $(z,0)$ has no orthogonal in $\overline{\Gamma}(R_1\times R_2)$,  a contradiction.\\
	$\Leftarrow)$ We have $Z(R_1\times R_2)= (R_1\setminus Z(R_1)\times Z(R_2))\cup   (Z(R_1) \times R_2\setminus Z(R_2))\cup (Z(R_1) \times Z(R_2)) $. Let $(a,b)\in Z(R_1\times R_2)$. If $(a,b)\in R_1\setminus Z(R_1)\times Z(R_2)$, then $(a,b)\perp (0,1)$. If $(a,b)\in Z(R_1)\times R_2\setminus Z(R_2)$, then $(a,b)\perp (c,0)$ with $c\in a^\perp$. If $(a,b)\in Z(R_1)\times Z(R_2)$, then $(a,b)\perp (c,1)$ with $c\in a^\perp$.
	\cqfd

 For the case where $R_2$ is a non reduced  ring such that  $|Z(R_2)|\geq 3$, we give the following theorem. 
\begin{thm}\label{products3}
	Let $R_1$ be a ring  and $R_2$ be a non reduced  ring such that  $|Z(R_2)|\geq 3$. Then, $\overline{\Gamma}(R_1\times R_2)$ is complemented if and only if $\overline{\Gamma}(R_2)$  and  $\overline{\Gamma}(R_1)$ are both complemented with $R_1$ is reduced.
\end{thm}
\pr 
$\Rightarrow)$ If $R_1$ is not reduced, then $\nil(R_1\times R_2)\geq 3$ since $R_2$ is not reduced. Then, $\overline{\Gamma}(R_1\times R_2)$ is not  complemented, by Theorem \ref{thmNil(R)}, since $|Z(R_1\times R_2)|\geq 4$, a contradiction.\\ Now, assume that $\overline{\Gamma}(R_2)$ is not complemented.  Then, there exists $z\in Z(R_2)^*$ such that for every $x\in Z(R_2)^*$, $x-z$ is a part of a triangle.
 Let $(a,b)\in Z(R_1\times R_2)^*$ such that $(a,b)$ is adjacent to $(0,z)$. Then, $(a,b)^n(0,z)^m=(0,0)$ for some $n,m\in \N^*$   with $(a,b)^n\neq (0,0)$ and $(0,z)^m\neq (0,0)$. Thus, $b^nz^m=0$ and so we have two cases to discuss:\\
 \textbf{Case 1.} If $b^n\neq 0$, then $b$ is adjacent to $z$ and so there exist a vertex $x$ adjacent to both $z$ and $b$. Thus, $(0,x)$ is adjacent to both $(0,z)$ and $(a,b)$.\\
 \textbf{Case 2.} If $b^n=0$, then $a^n\neq 0$ and so for every vertex $y\in Z(R_2)^*$ such that $y$ is adjacent to $z$, $(0,y)$ is adjacent to both $(a,b)$ and $(0,z)$. A contradiction since   $\overline{\Gamma}(R_1\times R_2)$ is complemented. Similarly,  we can prove that  $\overline{\Gamma}(R_1)$ is complemented, since if $\overline{\Gamma}(R_1)$ is  not complemented, $|Z(R_1)|\geq 3$.\\ 
 $\Leftarrow)$ We have $Z(R_1\times R_2)=(Z(R_1)\times Z(R_2))\cup (R_1\setminus Z(R_1)\times Z(R_2))\cup (Z(R_1)\times R_2\setminus Z(R_2))$. Let $(a,b)\in Z(R_1\times R_2)^*$. If $(a,b)\in R_1\setminus Z(R_1)\times Z(R_2)\setminus \nil(R_2)$, then $(a,b)\perp (0,c)$ with $c\in b^\perp$. If $(a,b)\in Z(R_1)^*\times R_2\setminus Z(R_2)$, then $(a,b)\perp (c,0)$ with $c\in a^\perp$. If $(a,b)\in Z(R_1)^*\times Z(R_2)\setminus \nil(R_2)$, then $(a,b)\perp (c_1, c_2)$ with $c_1\in a^\perp$ and $c_2\in b^\perp$. If $(a,b)\in R_1\setminus Z(R_1)\times \nil(R_1)^*$, then $(a,b)\perp (0,c)$  with $c\in R_2\setminus Z(R_2)$. If $(a,b)\in Z(R_1)^*\times \nil(R_2)^*$, then $(a,b)\perp (c_1,c_2)$ with $c_1\in a^{\perp}$ and $c_2\in R_2\setminus Z(R_2)$. If $(a,b)\in \{0\}\times \nil(R_2)^*$, then $(a,b)\perp  (c_1,c_2)$ with $c_1\in R_1\setminus Z(R_1)$ and $c_2\in b^\perp$.  
\cqfd 	
	
%	\begin{exm}
%		\begin{enumerate}
%			\item $\overline{\Gamma}(F_1\times F_2)$ and $\overline{\Gamma}(F_1\times \Z_4)$   are    complemented, where  $F_1$ and  $F_2$ are two fields.       
%			\item  $\overline{\Gamma}(\Z_6\times \Z_4)$,
%			$\overline{\Gamma}(\Z_6\times \Z_2)$,
%			$\overline{\Gamma}(\Z_{12}\times \Z_3)$ and 
%			$\overline{\Gamma}(\Z_6\times \Z_6)$  are complemented.
%		\end{enumerate}
%	\end{exm}

	Recall  that  the idealization of an $R-$module $M$, denoted by $R(+)M$, is the commutative ring $R\times M$ with the following addition and multiplication: $(a,n)+(b,m)=(a+b,n+m)$ and $(a,n)(b,m)=(ab,am+bn)$  for every   $(a,n),(b,m)\in  R(+)M$, \cite{J}. In  the following result we study when $\overline{\Gamma}(R(+)M)$ is complemented. Notice that, if $|M|\geq 4$, then $|Z(R(+)M)|\geq 4$ and $|\nil(R(+)M)|\geq 3$ and so $\overline{\Gamma}(R(+)M)$ is not complemented, by Theorem \ref{thmNil(R)}. If $M\cong \Z_3$, then $|\nil(R(+)\Z_3)|\geq 3$  and so $\overline{\Gamma}(R(+)\Z_3)$ is complemented if and only if $R$ is integral domain and $\Z_3$ is a torsion free $R$-module. Namely,  $\overline{\Gamma}(R(+)\Z_3)$ is an edge. Then, only the case where $M\cong \Z_2$ is of interest.

	\begin{pro}\label{idealization}
		Let $R$ be a  ring which is not an integral domain such that for every $x\in Z(R)\cap Z(\Z_2)$,  $x^\perp \setminus Z(\Z_2) \neq \emptyset$. Then,  $\overline{\Gamma}(R(+)\Z_2)$ is complemented if and only if $R$ is  reduced and  $\overline{\Gamma}(R)$ is complemented.
	\end{pro}
	\pr
	$\Leftarrow)$  We have  $Z(R(+)\Z_2)=Z(R)\cup Z(\Z_2)(+)\Z_2=\{(a,\bar{0}), (a,\bar{1})|\ a\in Z(R)\cup Z(\Z_2)\}$. Let $a\in Z(R)\cup Z(\Z_2)$, then we have the following three cases:\\
	\textbf{Case 1.} If $a\in Z(\Z_2)\setminus Z(R)$. Then,  $(a,\bar{0})\perp (0,\bar{1})$ and $(a,\bar{1})\perp (0,\bar{1})$.\\
	\textbf{Case 2.} If $a\in Z(R)\setminus Z(\Z_2)$. Then,  since $\overline{\Gamma}(R)$ is complemented, $(a,\bar{0})\perp (x,\bar{0})$ with $x\in  a^{\perp}$ and we have either $(a,\bar{1})\perp (y,\bar{0})$,  with $y\in a^{\perp}\cap Z(\Z_2)$, or $(a,\bar{1})\perp (y,\bar{1})$, with $y\in a^{\perp}\setminus Z(\Z_2)$.\\
	\textbf{Case 3.} If $a\in Z(R)\cap Z(\Z_2)$. Then,  there exists $x\in a^\perp \setminus Z(\Z_2)$ with $(a,\bar{0})\perp (x,\bar{0})$ and $(a,\bar{1})\perp (x,\bar{0})$.\\
	$\Rightarrow)$ Assume that $\overline{\Gamma}(R(+)\Z_2)$ is complemented. Then,  $|\nil(R(+)\Z_2)|=2$ (by Theorem  \ref{thmNil(R)}), which implies that $R$ is reduced. Now, let us prove that $\overline{\Gamma}(R)$ is complemented. Let $a\in Z(R)^*$. We have two cases:\\
	 \textbf{Case 1.} If $a\in Z(R)\cap Z(\Z_2)$, then $a$ has an orthogonal, by the hypotheses.\\
	  \textbf{Case 2.} If $a\in Z(R)\setminus Z(\Z_2)$, then   $(a,\bar{0})\in Z(R(+)\Z_2)^*$ and so $(a,\bar{0})$ has an orthogonal in $\overline{\Gamma}(R(+)\Z_2)$. Since the vertices adjacent to $(a,\bar{0})$ are  of the form $(b,\bar{1})$ or $(b,\bar{0})$ with $b\in Z(R)^*$, $(a,\bar{0})\perp (c,\bar{0})$ or $(a,\bar{0})\perp (c,\bar{1})$ for some $c\in a^\perp$. Hence,  $a$ has $c$ as an orthogonal.
	\cqfd

\end{document}